\documentclass[12pt, reqno]{amsart}
\usepackage{ifthen,amsfonts,amsmath,amssymb,epic,eepic}
\usepackage{epsfig,color}

\makeatother

\theoremstyle{definition}

\def\fnum{equation}
\newtheorem{Thm}[\fnum]{Theorem}
\newtheorem{Cor}[\fnum]{Corollary}

\newtheorem{Con}[\fnum]{Conjecture}

\newcommand{\nn}{{\bf{n}}}

\def\RR{{\bold R}}

\def\CC{{\bold C }}

\newcommand{\Area}{{\text {Area}}}

\newcommand{\K}{{\text{K}}}

\newcommand{\cP}{{\mathcal{P}}}

\newcommand{\eqr}[1]{(\ref{#1})}

\begin{document}

\title[What are the shapes of embedded minimal surfaces and why?]
{What are the shapes of embedded minimal surfaces and why?}

\author{Tobias H. Colding}%
\address{MIT and Courant Institute of Mathematical Sciences}
\author{William P. Minicozzi II}%
\address{Johns Hopkins University}
\thanks{The authors were partially supported by NSF Grants DMS
0104453 and DMS 0405695}


\email{colding@math.mit.edu and minicozz@math.jhu.edu}

\renewcommand{\abstractname}{Summary}
\begin{abstract}
Minimal surfaces with
uniform curvature (or area) bounds have been well understood and the regularity
theory is complete, yet essentially nothing was known without such
bounds.  We discuss here the theory of embedded 
(i.e., without self-intersections) minimal surfaces in 
Euclidean space $\RR^3$ 
without a priori bounds.  The study is divided into three cases, depending 
on the topology of the surface.  Case one is where the surface is a disk, 
in case two the surface is a planar domain (genus zero), and the third case is 
that of finite (non-zero) genus.  The complete 
understanding of the disk case is applied in both cases two and three.  

As we will see, the helicoid, which is a double
spiral staircase, is the most important example of an embedded minimal disk.  
In 
fact, we will see that every such disk is 
either a graph of a function or part of a
double spiral staircase.  The helicoid was discovered to be a minimal surface 
by Meusnier in 1776.  

For planar domains the fundamental examples are the catenoid, 
also discovered by Meusnier in 1776, and the 
Riemann examples discovered by Riemann in the beginning of the 
1860s\footnotemark.  Finally, 
for general fixed genus an important example is 
the recent example by Hoffman-Weber-Wolf of a genus one helicoid. 

In the last section we discuss why embedded minimal surfaces 
are automatically proper.  This was known as the Calabi-Yau conjectures 
for embedded surfaces.  For immersed surfaces there are 
counter-examples by Jorge-Xavier and Nadirashvili.
\end{abstract}
\footnotetext{Riemann worked on minimal surfaces in the period 1860-1861.  
He died in 1866.  
The Riemann example was published post-mortem in 1867 in an article edited
by Poggendorf.}

\maketitle


\numberwithin{equation}{section}


\setcounter{section}{0}

\setcounter{footnote}{1}
\subsection{Shape of things that are in equilibrium}
\text{  }

\vskip2mm
\centerline{What are the possible shapes of natural objects in equilibrium 
and why?}
  
\vskip2mm
\noindent
When a closed wire or a frame 
is dipped into a soap solution and is raised up from 
the solution, the surface spanning the wire is a soap film.  The soap film 
is in a state of equilibrium.     
What are the 
possible shapes of soap films and why?  Or  
why is DNA like a double spiral 
staircase?  ``What..?'' and ``why..?'' are fundamental questions which, 
when answered,  
help us understand the world we live in.  The answer to 
any question about the shape of natural objects is bound to 
involve mathematics.    

Soap films, soap bubbles, 
and surface tension were extensively studied by the Belgian physicist 
and inventor (of the stroboscope) Plateau 
in the first half of the 
nineteenth century.  At least since his studies, it has been known 
that the right  
mathematical model for a soap film 
is a minimal surface\footnote{The field of minimal surfaces  
dates back to the publication in 1762 of 
Lagrange's famous memoir ``Essai d'une nouvelle 
m\'ethode pour d\'eterminer les maxima et les minima des formules 
int\'egrales ind\'efinies''.  Euler had already, in a paper published in 
1744, discussed minimizing properties of the surface now known as the 
catenoid, but he only considered variations within a 
certain class of surfaces.} 
 -- the soap film is in a state of 
minimum energy when it is covering the least possible amount of area.  

We will discuss here the answer to the question: ``What are the possible 
shapes of embedded minimal surfaces in $\RR^3$ and why?''

\subsection{Critical points, minimal surfaces}

Let $\Sigma\subset \RR^3$ be a \underline{smooth} 
orientable surface (possibly with
boundary) with unit normal $\nn_{\Sigma}$.  Given
a function $\phi$ in the space $C^{\infty}_0(\Sigma)$ 
of infinitely differentiable (i.e., smooth), compactly supported 
functions on $\Sigma$, consider the one-parameter variation
\begin{equation}
\Sigma_{t,\phi}=\{x+t\,\phi (x)\,\nn_{\Sigma}(x) | x\in \Sigma\}\, .
\end{equation}
The so-called first variation formula of area is the equation (integration is
with respect to the area of $\Sigma$)
\begin{equation}  \label{e:frstvar}
\left.\frac{d}{dt} \right|_{t=0}\Area (\Sigma_{t,\phi})
=\int_{\Sigma}\phi\,H\, ,
\end{equation}
where $H$ is the mean curvature of $\Sigma$ and the mean curvature 
is the sum of the principal curvatures $\kappa_1$, $\kappa_2$.  
(When $\Sigma$ is non-compact, 
$\Sigma_{t,\phi}$ in \eqr{e:frstvar} is replaced by 
$\Gamma_{t,\phi}$, where 
$\Gamma$ is any compact set containing the support of $\phi$.)  
The surface $\Sigma$ is said to be a {\it minimal} surface (or just minimal) if
\begin{equation}
\left.\frac{d}{dt} \right|_{t=0}\Area (\Sigma_{t,\phi})=0
\,\,\,\,\,\,\,\,\,\,\,\text{ for all } \phi\in C^{\infty}_0(\Sigma)
\end{equation}
 or, equivalently by \eqr{e:frstvar}, if the
mean curvature $H$ is identically 
zero.  Thus $\Sigma$ is minimal if and only if
it is a critical point for the area functional.  Moreover, when $\Sigma$ 
is minimal, $\kappa_1=-\kappa_2$ (since $H=\kappa_1+\kappa_2=0$) and the 
Gaussian curvature $\K_{\Sigma}=\kappa_1\,\kappa_2$ is non-positive.

\subsection{Minimizers and stable minimal surfaces}

Since a critical 
point is not necessarily a minimum the term 
``minimal'' is
misleading, but it is time-honored.  The equation for a critical 
point is also sometimes called the Euler-Lagrange equation.   
A computation shows that if $\Sigma$ is minimal, then
\begin{equation}
\left. \frac{d^2}{dt^2} \right|_{t=0}\Area (\Sigma_{t,\phi})
=-\int_{\Sigma}\phi\,L_{\Sigma}\phi\, ,
\,\,\,\,\,\,\,\,\,\,\,\text{ where }
L_{\Sigma}\phi=\Delta_{\Sigma}\phi+|A|^2\phi
\end{equation}
is the second variational (or Jacobi) operator.
Here $\Delta_{\Sigma}$ is the Laplacian on $\Sigma$ and
$A$ is the second fundamental form of $\Sigma$.  So $A$ is the covariant 
derivative of the unit normal of $\Sigma$ and 
$|A|^2=\kappa_1^2+\kappa_2^2=-2\,\kappa_1\,\kappa_2=-2\K_{\Sigma}$,
 where $\kappa_1,\,\kappa_2$ are the principal curvatures 
(recall that since $\Sigma$ is minimal $\kappa_1=-\kappa_2$).
A minimal surface $\Sigma$ is said to be
stable if
\begin{equation}
\left. \frac{d^2}{dt^2} \right|_{t=0}\Area (\Sigma_{t,\phi})\geq 0
\,\,\,\,\,\,\,\,\,\,\,\text{ for all } \phi\in C^{\infty}_0(\Sigma)\, .
\end{equation}
One can show that a minimal graph is stable and, more generally, so is a
multi-valued minimal graph (see below for the precise definition).

Throughout, let $x_1 , x_2 , x_3$ be the standard coordinates on $\RR^3$.
For $y \in \Sigma \subset \RR^3$ and $s > 0$, the
extrinsic balls are
$B_s(y)=\{x\in \RR^3||x-y|<s\}$.

\subsection{Embedded = without self-intersections}

Our surfaces will all be without self-intersections, i.e., they will be 
embedded.  By embedded we 
mean a smooth injective immersion from an abstract surface into $\RR^3$.  

\subsection{Topology of surfaces}

The classification of minimal surfaces will essentially only depend on 
the topology of the surface and on whether or not the surface has a point 
where the curvature is large.

Compact orientable surfaces without boundaries are classified by their genus, a
nonnegative
integer. $\text{Genus}=0$ corresponds to a sphere, 
$\text{genus}=1$ to the torus, 
a model of which
is  the surface
of an object formed by attaching a ``suitcase handle'' to a sphere. 
A surface of
$\text{genus}=k$ is modelled by
the surface of a sphere to which $k$-handles have been attached. A compact
orientable surface with
boundary is one formed by taking  one of these surfaces and 
removing a number of disjoint 
disks.  The genus
of the surface with boundary is the genus of the original object, 
and the boundary
corresponds to the edges
of the surface created by disk removal.    
In particular, a surface with genus $0$ and non-empty boundary is a planar 
domain, i.e., it can be obtained from the disk in the plane 
by removing a number of 
disjoint sub-disks.  This is because it can be obtained from the sphere 
by removing a number of disks and after removing the first disk 
from the sphere, we have a disk in a plane.   Sometimes we will talk 
about surfaces that are simply connected.  By this we will mean that 
every loop on the surface can be shrunk (without leaving the surface) 
to a point curve.  One can easily see that the only simply connected 
surfaces are the disk and the sphere.  

\section{Disks}

There are two local models for embedded minimal disks.  One model is the plane
(or, more generally, a minimal graph) and the other is a piece of a 
helicoid.

\subsection{Minimal graphs and the helicoid}  \label{s:shelicoid}

The derivation of the equation for a minimal graph goes back 
to Lagrange's 1762 memoir.  (Note that if $\Omega$ is a simply 
connected domain in $\RR^2$ and 
$u$ is a real valued function, the graph of $u$, i.e., 
the set $\{(x_1,x_2,u(x_1,x_2))\,|\,
(x_1,x_2)\in \Omega\}$, is a disk.)  This gives a large class of 
embedded minimal disks where 
the function is defined over a proper subset of $\RR^2$, however by 
a classical theorem of Bernstein from 1916 entire (i.e., where
$\Omega=\RR^2$) minimal graphs are planes.  

The second model comes from the helicoid which was discovered by 
Meusnier in 1776\footnote{Meusnier had been a student of Monge.  He also 
discovered that the surface now known as the catenoid is minimal 
in the sense of Lagrange, and he was  
the first to characterize a minimal surface as a surface with 
vanishing mean curvature.  Unlike the helicoid, the catenoid is 
not topologically a plane but rather a cylinder.  
(The catenoid will be explained later; see \eqr{e:cat}.)}.
The 
helicoid is a ``double spiral staircase'' given by sweeping out a 
horizontal line rotating at a constant rate as it moves up a vertical 
axis at a constant rate.  Each half-line traces out a spiral 
staircase and together the two half-lines trace out (up to scaling) 
the double 
spiral staircase
\begin{equation}  \label{e:helicoid}
(s\cos t,s\sin t,t)\, ,\,\,\,\,\,\text{ where }s,\,t\in \RR\, .
\end{equation}

Anyone who has climbed a spiral staircase knows that the stairs 
become steep in the center.  The tangent plane to the helicoid at 
a point on the vertical axis is a vertical plane; thus the helicoid 
is not a graph over the horizontal plane.  In fact, as we saw earlier any 
minimal surface has non-positive curvature, for the 
helicoid the curvature is most negative along the axis and converges 
asymptotically to zero as one moves away from the axis. This 
corresponds to that as one moves away from the axis larger and 
larger pieces of the helicoid are graphs.  

For the our results about embedded minimal disks 
(see Sub-section \ref{s:sdisks} below) it will be 
important to understand a sequence of helicoids obtained from a single 
helicoid by rescaling as follows:  

\parbox{6in}{Consider the sequence 
$\Sigma_i = a_i \, \Sigma$ of 
rescaled helicoids where $a_i \to 0$. (That is, rescale $\RR^3$ by $a_i$, so 
points that used to be distance $d$ apart will in the rescaled $\RR^3$ be 
distance $a_i\,d$ apart.)  The curvatures of this 
sequence of rescaled helicoids 
are blowing up (i.e., the curvatures go to infinity) 
along the vertical axis. The sequence
converges (away from the vertical axis) to a foliation by flat
parallel planes. The singular set (the axis) then consists
of removable singularities.}

\subsection{Multi-valued graphs, spiral staircases, 
double spiral staircases} 
To be able to give a 
precise meaning to the statement 
that the helicoid is a double spiral staircase we will need
the notion of a multi-valued graph, each
staircase will be a multi-valued graph.  
Intuitively, a multi-valued graph is a 
surface covering an annulus, such that over a neighborhood of 
each point of the annulus, 
the surface consists of $N$ graphs. To make this notion 
precise, let $D_r$ be the disk in the plane 
centered at the origin and of radius 
$r$ and let $\cP$ be the universal cover of the
punctured plane $\CC\setminus \{0\}$ with global polar coordinates
$(\rho, \theta)$ so $\rho>0$ and $\theta\in \RR$.  An {\it $N$-valued
graph} on the annulus $D_s\setminus D_r$ is a
single valued graph of a function $u$
over $\{(\rho,\theta)\,|\,r< \rho\leq s\,
,\, |\theta|\leq N\,\pi\}$.  For working purposes, we generally think of the 
intuitive picture of a multi-sheeted surface in 
$\RR^3$, and we identify the single-valued graph over 
the universal cover with its multi-valued image in $\RR^3$.  

The multi-valued graphs that
we will consider will all be
embedded, which corresponds to a non-vanishing separation
between the sheets (or the floors).  
If $\Sigma$ is the helicoid, then
$\Sigma\setminus \{x_3-\text{axis}\}=\Sigma_1\cup \Sigma_2$,
 where $\Sigma_1$, $\Sigma_2$ are $\infty$-valued graphs on 
$\CC\setminus \{0\}$.
$\Sigma_1$ is the graph of the function $u_1(\rho,\theta)=\theta$
and $\Sigma_2$ is the
graph of the function $u_2(\rho,\theta)=\theta+\pi$.  ($\Sigma_1$ is the 
subset where $s>0$ in \eqr{e:helicoid} and $\Sigma_2$ the subset 
where $s<0$.)  In either
case the separation between the sheets is constant, equal to $2\,\pi$.  
A {\it multi-valued minimal graph} is a multi-valued graph of a
function $u$ satisfying the minimal surface equation.

\subsection{Structure of embedded minimal disks} \label{s:sdisks}
All of our results, for disks as well as for other topological types, 
require only a piece of a minimal surface.  In particular, the surfaces may 
well have boundaries.  This is a major point and makes the results 
particularly useful.
  
The following is the main structure theorem for embedded minimal 
disks\footnote{See \cite{CM1}--\cite{CM4} for the precise statements, 
as well as proofs, 
of the results of this section.}:

\begin{Thm} \label{t:t0.1}
Any embedded minimal disk in $\RR^3$ 
is \underline{either} a graph of a function \underline{or} part of a 
double spiral staircase.  In particular, if for some point the curvature 
is sufficiently large, then the surface is part of a double spiral staircase 
(it can be approximated by a piece of a rescaled helicoid).  On the 
other hand, if the curvature is below a certain threshold everywhere, 
then the surface is a 
graph of a function.
\end{Thm}

As a consequence of this structure theorem we get the 
following compactness result:

\begin{Cor}  \label{c:c1}
A sequence of embedded 
minimal disks with curvatures blowing up (i.e., going to 
infinity\footnote{Recall that for a minimal surface in $\RR^3$ the curvature 
$K=-\frac{1}{2}|A|^2$ is non-positive; so by that the curvatures 
of a sequence is going to infinity we mean that $K\to -\infty$ 
or equivalently $|A|^2\to \infty$.}) at a point  
mimics the behavior of a sequence of rescaled helicoids with curvature going 
to infinity; see the discussion of rescaled helicoids 
at the end of Sub-section \ref{s:shelicoid}.
\end{Cor}

\subsection{A consequence for sequences that are ULSC}
Sequences of planar domains which are not simply connected are,
after passing to a subsequence, naturally divided into two
separate cases depending on whether or not the topology is
concentrating at points.  To distinguish between these cases, we
will say that a sequence of surfaces $\Sigma_i^2\subset \RR^3$ is
{\it{uniformly locally simply connected}} (or ULSC)  if for each
$x\in\RR^3$, there exists a constant $r_0 > 0$
(depending on $x$) so that for all $r \leq
r_0$, and every surface
$\Sigma_i$
\begin{equation}    \label{e:ulsc2}
 {\text{each connected component of }} B_{r}(x) \cap \Sigma_i {\text{ is
 a disk.}}
\end{equation}
  For instance, a sequence of rescaled catenoids
   where the necks shrink to zero is not ULSC, whereas a
sequence of rescaled helicoids is.

The catenoid is the minimal surface in $\RR^3$
given by rotating the curve $s\to (\cosh s, s)$ around the $x_3$-axis, 
i.e., the surface
\begin{equation}    \label{e:cat}
 (\cosh s\, \cos t,\cosh s\, \sin t,s)  \text{ where }s,\,t\in \RR\,
 .
\end{equation}

\vskip2mm Applying the above structure theorem for disks to ULSC 
sequences gives that there are only two
{\underline{local models}} for such surfaces.  That is, 
locally in a ball in $\RR^3$, one of following
holds:
\begin{itemize}
\item  The curvatures are bounded and the surfaces are locally
{\underline{graphs}} over a plane.
\item
The curvatures blow up  and the surfaces are locally
{\underline{double spiral staircases}}.
\end{itemize}
Both of these cases are illustrated by taking a sequence of
rescalings of the helicoid; the first case occurs away from the
axis, while the second case occurs on the axis. 
If we take a sequence $\Sigma_i = a_i \, \Sigma$ of rescaled
helicoids where $a_i \to 0$, then the curvature blows up along the
vertical axis but is bounded away from this axis.  Thus, we get
that
\begin{itemize} \item The intersection of the rescaled helicoids
with a ball {\underline{away from}} the vertical axis gives a
collection of graphs over the plane $\{ x_3 = 0 \}$.
\item The intersection of the
rescaled helicoids with a ball {\underline{centered on}} the
vertical axis gives a double spiral staircase.
\end{itemize}

\subsection{Two key ideas behind the proof of the structure 
theorem for disks}
The first of these key ideas 
says that if the curvature of such a disk $\Sigma$
is large at some point $x\in \Sigma$, then near $x$
a multi-valued graph forms (in $\Sigma$), and this extends
(in $\Sigma$) almost all the way to the boundary of 
$\Sigma$\footnote{Recall that our results require only 
that we have a piece of a 
minimal surface and thus it may have boundary.}.
Moreover, the inner radius, $r_x$, of the annulus where the
multi-valued graph is defined is inversely proportional to $|A|(x)$,
and the initial separation between the sheets is bounded by a constant times
the inner radius.

An important ingredient in the proof of Theorem \ref{t:t0.1} is
that general embedded minimal disks with
large curvature at some interior point can be built out of
$N$-valued graphs.  In other words, any
embedded minimal disk can be divided into pieces each of which is
an $N$-valued graph.  Thus the disk itself should be thought of as
being obtained by stacking these pieces (graphs) on top of
each other\footnote{The parallel to the helicoid is striking.  Half of 
the helicoid, i.e., $(s\cos t,s\sin t, t)$, where $s>0$ and $t\in \RR$, 
can be obtained by stacking the $N$-valued graphs, 
$(s\cos (k\,N\,2\pi+t), s\sin (k\,N\,2\pi+t), k\,N\,2\pi+t)$, 
where $s>0$, $N\,2\pi>t\geq 0$, 
and $k$ is an integer, 
on the top of each other.}.

The second key result (Theorem \ref{t:t2}) is a curvature estimate for
embedded minimal disks in a half-space.  As a corollary 
(Corollary \ref{c:conecor} below) of this theorem,
we get that
the set of points in an embedded minimal disk where the curvature is
large lies within
a cone, and thus the multi-valued
graphs, whose existence was discussed above, will all start off within
this cone.

The curvature
estimate for disks in a half-space is the following:

\begin{Thm}  \label{t:t2}
There exists $\epsilon>0$ such that for all $r_0>0$, if
$\Sigma \subset B_{2r_0} \cap \{x_3>0\}
\subset \RR^3$ is an
embedded minimal
disk with $\partial \Sigma\subset \partial B_{2 r_0}$,
then for all components $\Sigma'$ of
$B_{r_0} \cap \Sigma$ which intersect $B_{\epsilon r_0}$
\begin{equation}  \label{e:graph}
\sup_{x\in\Sigma'} |A_{\Sigma}(x)|^2
\leq r_0^{-2} \, .
\end{equation}
\end{Thm}

Theorem \ref{t:t2} is an interior estimate where the curvature bound,
\eqr{e:graph},
is on the ball $B_{r_0}$ of one half of the radius of the ball $B_{2r_0}$
containing $\Sigma$.
This is just like a gradient estimate for a harmonic function where the
gradient bound is on one half of
the ball where the function is defined.

Using the minimal surface equation and the fact that $\Sigma'$ has points
close to a plane, it is not hard to see that, for $\epsilon>0$
sufficiently small, \eqr{e:graph} is equivalent to the statement
that $\Sigma'$
is a graph over the plane $\{x_3=0\}$.

We will often refer to Theorem \ref{t:t2}
as \emph{the one-sided curvature estimate} (since $\Sigma$ is assumed
to lie on one side of a plane).
Note that the assumption in Theorem \ref{t:t2}
that $\Sigma$ is simply connected (i.e., that $\Sigma$ is a disk) is crucial,
as can be seen from the example of a rescaled catenoid. 
Rescaled catenoids converge (with multiplicity two) to
the flat plane.  Likewise, by considering
the universal cover of the catenoid, one sees that Theorem \ref{t:t2} requires
the disk to be embedded,
and not just immersed.

In the proof of Theorem \ref{t:t0.1}, the following (direct)
consequence of Theorem \ref{t:t2} (with the $2$-valued graph playing
the role of the plane $\{x_3=0\}$) is
needed. 

\begin{Cor}   \label{c:conecor}
If an embedded
minimal disk contains a $2$-valued graph over an annulus in a plane, 
then away from a cone with axis
orthogonal to the $2$-valued graph the disk
consists of multi-valued graphs over annuli in 
the same plane.  
\end{Cor}

By definition, if $\delta>0$, then 
the (convex) cone with vertex at the origin, cone angle $(\pi/2 -
\arctan \delta)$, and axis parallel to the $x_3$-axis is the set
\begin{equation}
\{x\in \RR^3 \mid
x_3^2 \geq
\delta^2\,(x_1^2+x_2^2) \}\, .
\end{equation}

\subsection{Uniqueness theorems}
Using the above structure theorem for disks 
Meeks-Rosenberg, \cite{MeRo}, proved 
that the plane and the helicoid are the only
complete properly embedded
simply-connected minimal surfaces in $\RR^3$ 
(the assumption of properness can in fact be removed by \cite{CM6}; 
see Section \ref{s:CY}).
Catalan had proved in 1842 that any complete ruled minimal surface is
either a plane or a helicoid.
A surface is said to be \emph{ruled} if it has the
parametrization
\begin{equation}
X(s,t)=\beta (t)+ s\,\delta (t),\qquad\text{ where }s,\,t\in \RR,
\end{equation}
and $\beta$ and $\delta$ are curves in $\RR^3$.  The curve $\beta (t)$
is called the \emph{directrix} of the surface, and a line having
$\delta (t)$ as direction vector is called a \emph{ruling}.  For the
helicoid in \eqr{e:helicoid}, the $x_3$-axis is a directrix,
and for each fixed $t$ the line $s\to (s\,\cos t,s\,\sin t,t)$ is a ruling.

For cylinders there is a corresponding uniqueness theorem.  Namely, 
combining \cite{Sc}, \cite{Cn} (see also \cite{CM7}), and \cite{CM6} 
one has that any complete 
embedded minimal cylinder in $\RR^3$ is a catenoid.  

Conjecturally similar uniqueness theorems should hold for other 
planar domains and surfaces of fixed (non-zero) genus; cf. \cite{MeP}, 
\cite{HWW}.

\section{Planar domains}

We describe next 
two main structure theorems for {\it{non-simply connected}}
embedded minimal planar domains.  (The precise statements of 
these results and their proofs can be found in \cite{CM5}.)

The first of these asserts that any such surface
{\underline{without}} small necks\footnote{By ``without small 
necks'' we mean that the intersection of the surface with all extrinsic 
balls with sufficiently small radii consists of simply connected components; 
cf. the notion of ULSC for sequences above.} 
can be obtained by gluing
together two oppositely-oriented double spiral staircases.  Note 
that when one glues two oppositely oriented double spiral staircases together, 
then one remains at the same level if one 
circles both axes.  

The second gives a ``pair of pants'' decomposition of any such surface
when there {\underline{are}} small necks, cutting the surface
along a collection of short curves. After
the cutting, we are left with graphical pieces that are defined
over a disk with either one or two sub-disks removed (a
topological disk with two sub-disks removed is called a ``pair of
pants'').

Both of these structures occur as different extremes in the
two-parameter family of minimal surfaces known as the Riemann
examples.

\subsection{The catenoid and the Riemann examples}

\vskip2mm
When the sequence is no longer ULSC, then there are other local
models for the surfaces.  The simplest example is a sequence of
rescaled catenoids.  

 A sequence of rescaled  catenoids converges with
multiplicity two to the flat plane.  The convergence is in the
$C^{\infty}$ topology except at $0$ where $|A|^2 \to \infty$. This
sequence of rescaled catenoids is not ULSC because the simple
closed geodesic on the catenoid -- i.e., the unit circle in the
$\{ x_3 = 0 \}$ plane -- is non-contractible and the rescalings
shrink it down to the origin.

One can get other types of curvature blow-up by considering the 
family   of embedded minimal planar domains known as the Riemann 
examples.  Modulo translations and rotations, this is a 
two-parameter family of periodic minimal surfaces, where the 
parameters can be thought of as the size of the necks and the 
angle  from one fundamental domain to the next. By choosing the 
two parameters appropriately, one can produce sequences of Riemann 
examples that illustrate both of the two structure theorems:

\begin{enumerate}

\item

If we take a sequence of Riemann examples where the neck size is 
fixed and the angles go to $\frac{\pi}{2}$, then the surfaces with 
angle near $\frac{\pi}{2}$ can be obtained by gluing together two  
oppositely-oriented double spiral staircases. Each double spiral 
staircase looks like a helicoid. This sequence of Riemann examples 
converges to a foliation by parallel planes.  The convergence is 
smooth away from the axes of the two helicoids (these two axes are 
the singular set where the curvature blows up). The sequence 
is ULSC since the size of the necks is fixed and  thus illustrates 
the first structure theorem, Corollary \ref{c:c2} below.

\item

If we take a sequence of examples where the neck sizes go to zero, 
then we get a sequence  that is {\it{not}} ULSC.  However, the 
surfaces can be cut along   short curves into collections of 
graphical pairs of pants.  The short curves converge to points and 
the graphical pieces converge to flat planes except at these 
points, illustrating the second structure theorem, Corollary  
\ref{c:c3} below. 
\end{enumerate}

\subsection{Structure of embedded planar domains}

We describe next (Theorems \ref{t:mst1} and \ref{t:mst2} below) 
the two main structure theorems for {\it{non-simply connected}}
embedded minimal planar domains.  Each of these theorems has a 
compactness theorem as a consequence.  

The first structure theorem deals with surfaces 
\underline{without} small necks:
 
\begin{Thm}  \label{t:mst1}
Any non-simply connected embedded minimal planar domain without small necks 
can be obtained from gluing together two oppositely oriented double 
spiral staircases. Moreover, if for some point the curvature 
is large, then the separation between the sheets of the double 
spiral staircases is small.  Note that since the two double spiral staircases 
are oppositely oriented, then one remains at the same level if one 
circles both axes.  
\end{Thm}

The following compactness result is a consequence:  

\begin{Cor}  \label{c:c2}
A ULSC (but not simply connected) sequence of embedded minimal surfaces 
with curvatures blowing up has 
a subsequence that converges smoothly to a
foliation by parallel planes away from two lines. The two lines are
disjoint and orthogonal to the leaves of the foliation, and the two
lines are precisely the points where the curvature is blowing up.

This is similar to the case of disks, except that we get two
singular curves for non-disks as opposed to just one singular
curve for disks.

Moreover, locally around each of the two lines the surfaces look 
like a helicoid around the axis and the orientation around the 
two axes are opposite.     
\end{Cor}

Despite the similarity of Corollary \ref{c:c2} 
 to the case of disks, it is worth noting that the results for
 disks do not alone give this result.  Namely, even though the
 ULSC sequence consists {\underline{locally}} of disks, the compactness
result for disks
 was in the {\underline{global}} case where the radii go to
 infinity.  One might wrongly assume that Corollary \ref{c:c2} could be proven
 using the results for disks and a blow-up argument. However, one can 
construct local
 examples that show the difficulty of such an
 argument.

The second structure theorem deals with surfaces 
\underline{with} small necks 
and gives a ``pair of pants'' decomposition:

\begin{Thm}  \label{t:mst2}
Any non-simply connected embedded minimal planar domain with a small 
neck can be cut 
along a collection of short curves. After
the cutting, we are left with graphical pieces that are defined
over a disk with either one or two sub-disks removed (a
topological disk with two sub-disks removed is called a pair of
pants).

Moreover, if for some point the curvature 
is large, then all the necks are very small.  
\end{Thm}

The following compactness result is a consequence:  

\begin{Cor}  \label{c:c3}
A sequence of embedded minimal planar domains that are not ULSC, but 
with curvatures blowing up, has a subsequence that 
converges to a lamination by flat parallel 
planes.  
\end{Cor}

\section{Finite genus}

\subsection{The genus one helicoid; structure results for general finite genus}

In a very recent paper Hoffman-Weber-Wolf, \cite{HWW}, have constructed a 
new complete embedded minimal surface in $\RR^3$.  They have 
shown that there exists a properly 
embedded minimal surface of genus one with a single end asymptotic 
to the end of the helicoid.  We will refer to this minimal surface $\Sigma$ 
as the genus one helicoid.  Under scalings the sequence of genus 
one surfaces $a_i\,\Sigma$ where $a_i\to 0$ 
converges to the foliation of flat parallel 
planes in $\RR^3$ just like a sequence of rescaled helicoids with 
curvatures blowing up. This is in fact a consequence of a  
general result that the theorems in the previous section,  
stated for planar domains, holds also
for sequences with fixed genus with minor changes; see \cite{CM5}.

\section{Embedded minimal surfaces are automatically proper} 
\label{s:CY}

Implicit in all of the results mentioned above was an assumption 
that the minimal surfaces were proper.  However, as we 
will see next, it turns out 
that embedded minimal surfaces are, in fact, automatically proper.  
This was the content of the Calabi-Yau conjectures which were proven 
to be true for embedded surfaces in \cite{CM6}.  

\subsection{What is proper?}

An immersed surface in $\RR^3$ is {\it proper} if the pre-image of any compact
subset of $\RR^3$ is compact in the surface.  For instance, a line 
is proper whereas a curve that spiral infinitely into a circle is not.   

\subsection{The Calabi-Yau conjectures; the statements and examples} 
The Calabi-Yau conjectures about surfaces date back to the 1960s.
Their original form  was given in 1965 where 
Calabi made the following two conjectures about minimal surfaces\footnote{S.S. 
Chern also promoted these conjectures at roughly the same time 
and they were revisited several times by S.T. Yau.}:

\begin{Con} \label{con:1}
``Prove that a complete  minimal surface in $\RR^3$ must be
unbounded.''
\end{Con}

Calabi continued:
``It is known that there are no compact minimal surfaces in 
$\RR^3$ (or of any simply connected complete Riemannian $3$-dimensional 
manifold
with sectional curvature $\leq 0$). A more ambitious conjecture
is'':

\begin{Con} \label{con:2}
``A complete [non-flat] minimal surface in $\RR^3$ has an
unbounded projection in every line.''
\end{Con}

 The \underline{immersed}  versions of these conjectures turned out to be
 false.  Namely,  Jorge and Xavier, \cite{JXa},  constructed
non-flat minimal immersions contained between two parallel planes
in 1980, giving a counter-example to the immersed version of the
more ambitious Conjecture \ref{con:2}.
 Another significant development came in 1996,  when
 Nadirashvili, \cite{Na}, constructed a complete
immersion of a minimal disk into the unit ball in $\RR^3$, showing
that Conjecture \ref{con:1} also failed for immersed surfaces; see 
\cite{LMaMo},
 for other topological types.

The main result in \cite{CM6} is an effective version of properness for disks,
giving a chord-arc bound\footnote{A chord-arc bound is a bound above and 
below for the ratio of intrinsic to extrinsic distances.}. Obviously, 
intrinsic distances are
larger than extrinsic distances, so the significance of a chord-arc 
bound is the reverse inequality, i.e., a bound on intrinsic
distances from above by extrinsic distances.  

 Given such a chord-arc bound, 
one has that as 
intrinsic distances go to infinity, so do extrinsic distances.
Thus as an immediate consequence:

\begin{Thm}     \label{t:2nd}
A complete embedded minimal disk in $\RR^3$ must be proper.
\end{Thm}

Theorem \ref{t:2nd} gives immediately that the first of
  Calabi's conjectures is true for
{\underline{embedded}} minimal disks. In particular,
Nadirashvili's examples cannot be embedded.
 
Another immediate consequence of the chord-arc bound together with
the one-sided curvature estimate (i.e., Theorem \ref{t:t2}) 
is a  version of that estimate
for intrinsic balls.

 As a corollary of this intrinsic one-sided curvature estimate we get that
the second, and more ambitious, of Calabi's conjectures is also
true for {\underline{embedded}} minimal disks.  In particular,
 Jorge-Xavier's examples cannot be embedded. The second Calabi conjecture 
(for embedded disks) is an immediate consequence of the 
following half-space
 theorem:

\begin{Thm}     \label{t:1}
The plane is the only complete embedded minimal disk in $\RR^3$
 in a half-space.
\end{Thm}

The results for disks imply
both of Calabi's conjectures and properness also for embedded
surfaces with finite topology. A surface $\Sigma$ is
said to have finite topology if it is homeomorphic to a closed
Riemann surface with a finite set of points removed or
``punctures''.  Each puncture corresponds to an end of $\Sigma$.

\vskip4mm
We thank C.H. Colding, L. Hesselholt, N. Hingston, and 
D. Hoffman for helpful comments.

\end{document}